\documentclass[12pt]{article}
\usepackage[utf8]{inputenc}
\usepackage{amssymb,amsmath,amsfonts,amsthm,amscd,latexsym,indentfirst,verbatim,xcolor}
\usepackage[T2A]{fontenc}
\usepackage{geometry}
\def\dbl{\lbrace\kern-3pt\lbrace}
\def\dbr{\rbrace\kern-3pt\rbrace}

\def\Imm{\textrm{Im}}
\def\End{\textrm{End}}
\def\Hom{\textrm{Hom}}

\def\As{\textrm{As}}
\def\Span{\textrm{Span}}
\def\id{\textrm{id}}

\def\tr{\textrm{tr}}

\begin{document}

\sloppy

\hfill{16W99, 17B63 (MSC2020)}

\begin{center}
{\Large An example of a simple double Lie algebra}

\smallskip

Vsevolod Gubarev
\end{center}

\begin{abstract}
We extend the correspondence between double Lie algebras and skew-symmetric
Rota---Baxter operators of weight 0 on the matrix algebra 
to the infinite-dimensional case. 
We give the first example of a simple double Lie algebra.

\end{abstract}

\section{Introduction}

In~2008~\cite{DoublePoisson},
M. Van den Bergh introduced the notion of a double Poisson algebra
developing noncommutative geometry.
For this, he followed the Kontsevich---Rosenberg principle
saying that a~structure on an associative algebra
has geometric meaning if it induces standard geometric structures
on its representation spaces.

Given a finitely generated associative algebra~$A$ and $n\in\mathbb{N}$,
consider the representation space $\mathrm{Rep}_n(A) = \Hom(A,M_n(F))$,
where $F$ denotes the ground field.
To equip~$A$ with a~structure such that
$\mathrm{Rep}_n(A)$ is a~Poisson variety for every~$n$,
M. Van den Bergh defined a double bracket
$\dbl \cdot,\cdot\dbr\colon A\otimes A\to A\otimes A$
satisfying the analogues of anti-commutativity,
Jacobi identity, and Leibniz rule.
An associative algebra equipped with such
a~double bracket is called a double Poisson algebra.
One of the crucial examples of such structure
is a double Poisson algebra defined on a quiver algebra.

Double Poisson algebras are deeply connected with
$H_0$-Poisson structures~\cite{Crawley-Boevey},
pre-Calabi---Yau algebras~\cite{IKV},
vertex algebras~\cite{Kac15}.

The notion of a double Lie algebra naturally appeared from the very definition
of double Poisson algebra, it is a vector space endowed with a double bracket
satisfying above mentioned anti-commutativity and Jacobi identity.
Every double Lie algebra structure defined on a vector space~$V$ 
can be uniquely extended to a~double Poisson algebra structure 
on the free associative algebra $\As\langle V\rangle$.
Thereby, A.~Odesskii, V.~Rubtsov, V.~Sokolov 
extended~\cite{DoublePoissonFree} linear and quadratic 
double Lie algebras defined on an $n$-dimensional vector space to
double Poisson algebras defined on the free $n$-generated associative algebra.

In~\cite{DoubleLie}, M. Goncharov and P. Kolesnikov proved
that there are no simple finite-dimensional double Lie algebras.
This problem was stated by V. Kac during the conference
``Lie and Jordan algebras, their representations and applications''
dedicated to Efim Zelmanov’s 60th birthday (Bento Gon\c{c}alves, Brasil, 2015).
After this work the natural question about constructing simple
infinite-dimensional double Lie algebras has arisen.

It is known that the structure of a~double Lie algebra on a finite-dimensional
vector space~$V$ is equivalent to a~skew-symmetric 
Rota---Baxter operator of weight~0 on
the matrix algebra~$M_n(F)$, where 
$n = \dim(V)$~\cite{DoubleLie,DoublePoissonFree,Schedler}.
Recall that a linear operator~$R$ defined on an algebra~$A$ is called
a~Rota---Baxter operator (RB-operator, for short) of weight~$\lambda$, if
$$
R(x)R(y) = R( R(x)y + xR(y) + \lambda xy )
$$
for all $x,y\in A$.
This notion for the first time appeared in the article~\cite{Tricomi} 
of F. Tricomi in 1951 and further was several times~\cite{Baxter,BelaDrin82} 
rediscovered, see the monograph~\cite{GuoMonograph}.
Let us mention the bijection~\cite{Aguiar00,Unital,Schedler} 
between RB-operators of weight~0 on the matrix algebra $M_n(F)$ 
and solutions of the associative Yang---Baxter
equation (AYBE) on $M_n(F)$~\cite{Aguiar01,Polishchuk,Zhelyabin}.

We generalize this correspondence between
double Lie algebras and skew-symmetric Rota---Baxter operators
for the infinite-dimensional case. We state such correspondence for
a countable-dimensional double Lie algebra~$V$ and
a Rota---Baxter operator acting from the space of matrices 
with finite numbers of nonzero elements to $\End(V)$
and satisfying some additional finiteness conditions.
This correspondence allows us to construct new double Lie algebras.
In particular, we show that the vector space $F[t]$ with the double bracket
$$
\dbl t^n,t^m\dbr = -\frac{(t^n\otimes t^m-t^m\otimes t^n)}{t\otimes 1-1\otimes t}
$$
is a simple double Lie algebra.
As far as we know it is the first example of a simple double Lie algebra.

In terms of RB-operators we interpret the amazing double 
Lie algebra of V.~Kac (see~\cite{DoubleLie})
whose definition is very close to the definition of the Yangian $Y(gl_N)$.

\section{Preliminaries}

\subsection{Rota---Baxter operators}

{\bf Definition 1}.
A~linear operator~$R$ defined on a (not necessary associative) algebra~$A$
is called a~Rota---Baxter operator (RB-operator, for short) 
of weight~$\lambda\in F$, if
\begin{equation} \label{RB}
R(x)R(y) = R( R(x)y + xR(y) + \lambda xy )
\end{equation}
holds for all $x,y\in A$.

{\bf Proposition 1}~\cite{Unital}.
Let $A$ be an algebra, let $R$ be an RB-operator of weight~$\lambda$ on~$A$,
and let $\psi$ be either automorphism or antiautomorphism of $A$.
Then the operator $R^{(\psi)} = \psi^{-1}R\psi$ is an RB-operator 
of weight~$\lambda$ on~$A$.

As an application of Proposition~1, we will use the conjugation with transpose
of an RB-operator defined on the matrix algebra.

The following definition has appeared by the name 
of relative Rota---Baxter operator
or $\mathcal{O}$-operator~\cite{O-operator} or as
generalized RB-operator in the case of zero weight~\cite{Uchino}. 
For simplicity, we also call it Rota---Baxter operator.

{\bf Definition 2}.
Let $A$ be an algebra and $I$ be an ideal of~$A$.
A~linear operator~$R\colon I\to A$ is called
a~Rota---Baxter operator of weight~$\lambda$, if
$$
R(i)R(j) = R( R(i)j + iR(j) + \lambda ij)
$$
holds for all $i,j\in I$.

When $I = A$, this definition coincides with Definition~1.

The next statement follows immediately.

{\bf Proposition 2}.
Let $A$ be an algebra and let~$J$ be an ideal of~$A$.
Given an RB-operator $P\colon J\to A$ of weight~$\lambda$ and an algebra~$B$,
the operator $Q = P\otimes \id_B$ is again an RB-operator
of the same weight~$\lambda$ from $I\otimes B$ to $A\otimes B$.

Henceforth, we consider only Rota---Baxter operators of weight~0.
It is well-known that, given an RB-operator~$R$ of weight~0
and $\alpha\in F$, the operator $\alpha R$
is again an RB-operator~$R$ of weight~0.

\subsection{Double Lie algebras}
Let $V$ be a linear space. Given $u\in V^{\otimes n}$ and $\sigma \in S_n$,
$u^\sigma$ denotes the permutation of tensor factors.
By a double bracket on $V$ we call a~linear map from $V\otimes V$ to $V\otimes V$.
Given an associative algebra~$A$, we consider 
the outer bimodule action of $A$ on $A\otimes A$:
$b(a\otimes a') c = (ba)\otimes (a'c)$.

{\bf Definition 3}~\cite{DoublePoisson}.
A double Poisson algebra is an associative algebra $A$ 
equipped with a~double bracket
satisfying the following identities for all $a,b,c\in A$
\begin{gather}
\dbl a,b\dbr =- \dbl b,a\dbr ^{(12)}, \label{antiCom} \\
\dbl a, \dbl b,c\dbr \dbr _L -\dbl b, \dbl a,c\dbr \dbr _R^{(12)}
 = \dbl \dbl a,b\dbr,c\dbr _L,  \label{Jacobi} \\
\dbl a,bc\dbr = \dbl a,b\dbr c + b\dbl a,c\dbr, \label{Leibniz}
\end{gather}
where
$\dbl a, b\otimes c \dbr _L = \dbl a,b \dbr \otimes c$,
$\dbl a, b\otimes c\dbr _R = (b\otimes \dbl a,c\dbr )^{(12)}$, and
$\dbl a\otimes b, c\dbr _L = (\dbl a,c\dbr \otimes b)^{(23)}$.

{\bf Definition 4}~\cite{DoublePoissonFree,Schedler,Kac15}.
A double Lie algebra is a linear space $V$ equipped with a~double bracket
satisfying the identities~\eqref{antiCom} and~\eqref{Jacobi}.

Due to~\cite{DoubleLie}, an {\it ideal}
of a~double Lie algebra $V$ is a subspace $I\subseteq V$ such that
$$
\dbl V,I\dbr + \dbl I,V\dbr \subseteq I\otimes V + V\otimes I.
$$
Given an ideal $I$ of a double Lie algebra~$V$, we have
a~natural structure of a double Lie algebra on the quotient space $V/I$, i.\,e.,
$\dbl x + I, y + I \dbr = \dbl x,y \dbr + I\otimes V + V\otimes I$.

Let us define homomorphisms of double Lie algebras as follows.
Let~$L$ and $L'$ be double Lie algebras and let $\varphi\colon L\to L'$
be a~linear map. Then~$\varphi$ is called a {\it homomorphism} from~$L$ to~$L'$ if
$$
(\varphi\otimes\varphi)(\dbl a,b\dbr) = \dbl \varphi(a),\varphi(b)\dbr
$$
holds for all $a,b\in L$.
Note that the kernel of any homomorphism from~$L$ is an ideal of~$L$.

{\bf Definition 5}~\cite{DoubleLie}.
A double Lie algebra $V$ is said to be simple if $\dbl V,V\dbr \neq (0)$
and there are no nonzero proper ideals in $V$.

\section{Finite-dimensional double Lie algebras}

Suppose that $V$ is a~finite-dimensional space.
In~\cite{DoubleLie}, it was shown that every double Lie algebra structure
$\dbl \cdot,\cdot\dbr$ on~$V$ is determined by a linear operator
$R\colon\End (V)\to \End(V)$, precisely,
\begin{equation}\label{eq:Bracket_via_RB}
 \dbl a,b\dbr  = \sum\limits_{i=1}^N e_i(a)\otimes R(e_i^*)(b)
 % = \sum\limits_{i=1}^N R^*(e_i)(a)\otimes e^*_i(b)
 , \quad a,b\in V,
\end{equation}
where $e_1,\dots, e_N$ is a linear basis of $\End(V)$, $e_1^*,\dots, e_N^*$
is the corresponding dual basis relative to the trace form.

A linear operator $P$ on~$\End(V)$ is called skew-symmetric if $P = -P^*$,
where $P^*$ is the conjugate operator on $\End(V)$ relative to the trace form.

{\bf Theorem~1}~\cite{DoubleLie}.
Let $V$~be a finite-dimensional vector space with 
a~double bracket $\dbl \cdot,\cdot \dbr$
determined by an operator $R\colon\End(V)\to \End(V)$ by~\eqref{eq:Bracket_via_RB}.
Then $V$ is a double Lie algebra if and only
if $R$ is a~skew-symmetric RB-operator of weight~0 on $\End(V)$.

{\bf Remark 1}.
Theorem~1 was stated in~\cite{Schedler} in terms of skew-symmetric
solutions of the associative Yang---Baxter equation (AYBE).
Since there is a one-to-one correspondence between
solutions of AYBE and Rota---Baxter operators 
of weight~0 on the matrix algebra~\cite{Unital},
Theorem~1 follows from~\cite{Schedler}. 
Actually, Theorem~1 was mentioned also in~\cite{DoublePoissonFree}.

Let us consider several examples of double Lie algebras 
and corresponding RB-operators.
We will use the linear basis $e_{ij}$, $1\leq i,j\leq \dim(V)$, of $\End(V)$.
So, we have $e_{ij}^* = e_{ji}$ relative to the trace form.

In the case of a~one-dimensional double Lie algebra~$L$, 
we have by~\eqref{antiCom} only zero double bracket.

{\bf Example 1}~\cite{DoubleLie,DoublePoissonFree,DoublePoisson}.
The space $F^2 = Fe_1 \oplus Fe_2$ equipped with a double product
$\dbl e_1,e_1\dbr = e_1\otimes e_2-e_2\otimes e_1$
(others are zero) is a double Lie algebra. The corresponding RB-operator on $M_2(F)$ is
$R_1 (e_{11}) = e_{21}$ and $R_1(e_{12}) = -e_{11}$ (others are zero).

{\bf Example 2}~\cite{DoubleLie,DoublePoissonFree,PichereauWeyer}.
The space $F^2$ with a double product
$\dbl e_1,e_2\dbr = e_1\otimes e_1 = - \dbl e_2,e_1\dbr$ 
is again a double Lie algebra.
The corresponding RB-operator on $M_2(F)$ is
$R_2 (e_{11}) = e_{12}$ and $R_2 (e_{21}) = -e_{11}$.

The RB-operators $R_1$ and $R_2$ are conjugate with the transpose of matrices,
i.\,e., $R_2 = R_1^{(T)}$, where $T$ denotes the transpose.
However, the algebraic properties of the double Lie algebras from
Examples~1 and~2 are quite different, see~\cite{DoubleLie}.

Note that all RB-operators (including skew-symmetric)
of weight~0 on $M_2(F)$ were classified by M.~Aguiar~\cite{Aguiar00-2} in 2000
and all skew-symmetric RB-operators of weight~0 on $M_3(\mathbb{C})$ 
were described by V.~V.~Sokolov~\cite{Sokolov} in 2013.

{\bf Example 3}.
Consider the restriction of the double bracket defined in~\cite[\S6.5]{DoublePoisson}
on the infinite-dimensional path algebra over a field $F$ arisen from the quiver~$Q$
with the vertex set $\{1,2\}$ and the edge set $\{e_1,e_2,a,a^*\}$, where
$a = (1,2)$ and $a^* = (2,1)$.
We put $L = \Span\{e_1,e_2,a,a^*\}$, and the double bracket on~$L$ equals
$$
\dbl a,a^*\dbr =  e_2\otimes e_1, \quad
\dbl a^*,a\dbr = -e_1\otimes e_2,
$$
all other double brackets are zero.
Let us identify $e_3 = a$ and $e_4 = a^*$.
By~\eqref{eq:Bracket_via_RB} we get the RB-operator $R$ on $M_4(F)$
defined as follows, $R(e_{32}) = e_{14}$, $R(e_{41}) = -e_{23}$.

\section{Infinite-dimensional double Lie algebras}

Consider a countable-dimensional double Lie algebra 
$\langle V,\dbl \cdot,\cdot\dbr\rangle$.
We fix a~linear basis $u_i$, $i\in\mathbb{N}$, of $V$.
Define $e_{ij}\in\End(V)$ by the formula $e_{ij}u_k = \delta_{jk}u_i$.
Let $\varphi\in\End(V)$, then we may write $\varphi = \sum\limits_{ij}a_{ij}e_{ij}$.
We identify $\varphi$ with an infinite matrix $[\varphi] = (a_{ij})_{i,j\geq0}$.
Since $\varphi\in\End(V)$ is well-defined,
there is only a~finite number of nonzero elements in every column 
of the matrix $[\varphi]$, i.\,e., $a_{ik} = 0$ for almost all $i$ when $k$~is fixed.

Let us define the subalgebra~$\End_f(V)$ of $\End(V)$ as follows,
$$
\End_f(V) = \{\varphi\in\End(V)\mid
 \mbox{ for every }i,\ [\varphi]_{ij} = 0
 \mbox{ for almost all }j\}.
$$

Introduce $I$ as an ideal in $\End_f(V)$ linearly spanned by matrix unities $e_{ij}$.

Let $\varphi\in\End_f(V) = \sum\limits_{i,j}a_{ij}e_{ij}$.
We define the symmetric non-degenerate bilinear trace form
$\langle \cdot,\cdot \rangle$ on $I\times \End_f(V)\cup \End_f(V)\times I$ as follows,
$$
\langle e_{kl},\varphi \rangle
 = \langle \varphi,e_{kl}\rangle
 = \tr(e_{kl}\varphi)
 = a_{lk}.
$$
Moreover, the form is associative, i.\,e., $\langle a,bc\rangle = \langle ab,c\rangle$,
where at least one of $a,b,c$  lies in~$I$ and others are from $\End_f(V)$.

Given a double bracket algebra~$\dbl \cdot,\cdot\dbr$ on a space~$V$, 
we may define a linear operator $R\colon I\to \End (V)$ by the formula
\begin{equation}\label{Bracket_via_RBInf}
\dbl a,b\dbr  = \sum\limits_{i,j\geq0} e_{ij}(a)\otimes R(e_{ji})(b), \quad a,b\in V.
\end{equation}
Conversely, given an operator $R\colon I\to \End (V)$, 
one can define a double bracket on~$V$ by the formula~\eqref{Bracket_via_RBInf}.
Note that the correspondence does not work if $R\colon \End(V)\to \End (V)$. 

Moreover, we define a conjugate operator $R^*\colon I\to \End (V)$ as follows,
\begin{equation}\label{Conjugate}
\dbl b,a\dbr^{(12)}  = \sum\limits_{i,j\geq0} e_{ij}(a)\otimes R^*(e_{ji})(b), 
\quad a,b\in V.
\end{equation}

Denote $R(e_{st}) = \sum\limits_{k,l}a_{kl}^{st}e_{kl}$. 
By~\eqref{Bracket_via_RBInf}, $a_{kl}^{st}$ equals the coefficient 
by~$u_t\otimes u_k$ of the double product $\dbl u_s,u_l\dbr$.
Analogously, put $R^*(e_{st}) = \sum\limits_{k,l}b_{kl}^{st}e_{kl}$. 
Then $b_{kl}^{st}$ equals the coefficient by~$u_k\otimes u_t$ 
of the double product $\dbl u_l,u_s\dbr$. Hence, 
$$
\langle R(e_{st}),e_{kl}\rangle 
 = \dbl u_s,u_k\dbr|_{u_t\otimes u_l} 
 = \langle R^*(e_{kl}),e_{st}\rangle.
$$
Generally we have 
\begin{equation}\label{InfConj}
\langle R(x),y\rangle = \langle x,R^*(y)\rangle,\quad x,y\in I.
\end{equation}

{\bf Remark 2}.
It is not clear how to introduce objects defined above in invariant manner.
For example, consider a linear basis $u_s$, $s\in\mathbb{N}$, 
of $V$ and a~linear map
$\psi\in\End_f(V)$ defined as follows,
$\psi(u_0) = u_0$ and $\psi(u_s) = 0$, $s>0$.
Let us consider the basis $w_s$, $s\in\mathbb{N}$, of $V$, where
$w_0 = u_0$ and $w_s = u_0 + u_s$, $s>0$.
Then $\psi(w_s) = w_0$ for all~$s$. Thus, $\psi\not\in\End_f(V)$. 
Hence, the change of the basis does not preserve the condition 
$R\colon I\to \End_f(V)$.

{\bf Theorem~2}.
Let $V$~be a countable-dimensional vector space with a fixed linear basis
$u_i$ and with a~double bracket
$\dbl \cdot,\cdot \dbr$ determined by a linear map $R\in \End_f(V)$
by~\eqref{Bracket_via_RBInf}.
Then $V$ is a double Lie algebra if and only if $R$ 
is a~skew-symmetric RB-operator
of weight~0 from~$I$ to $\End_f(V)$.

{\sc Proof}.
By~\eqref{Bracket_via_RBInf} and~\eqref{Conjugate},
the identity \eqref{antiCom} holds if and only if $R=-R^*$.

Define $F_{12}\in \End(V^{\otimes 3})$ by
$$
 F_{12}(a\otimes b\otimes c) = \dbl a, \dbl b,c\dbr \dbr_L
 = \sum\limits_{i,j} e_j(a)\otimes R(e_j^*)(e_i(b))\otimes R(e_i^*)(c),
 \quad a,b,c\in V.
$$
For $x,y\in I$, we compute applying associativity 
of the form $\langle\cdot,\cdot\rangle$
\begin{multline}\label{eq:F12}
 (\langle x,\cdot \rangle \otimes \langle y,\cdot \rangle \otimes \id)F_{12}
 = \sum\limits_{i,j} \langle x,e_j\rangle \langle y,R(e_j^*)e_i \rangle R(e_i^*) 
 \allowdisplaybreaks \\
 = \sum\limits_i\left\langle y,
 \sum\limits_j\langle x,e_j\rangle R(e^*_j)e_i\right\rangle R(e_i^*)
 = \sum\limits_i \langle y,R(x)e_i\rangle R(e_i^*) \\
 = \sum\limits_i \langle yR(x), e_i\rangle R(e_i^*)
 = R(yR(x)).
\end{multline}

Analogously, put
\begin{gather*}
F_{23}(a\otimes b\otimes c) = \dbl b, \dbl a,c\dbr \dbr _R^{(12)}
 = \sum\limits_{i,j} e_j(a)\otimes e_i(b) \otimes R(e_i^*)(R(e_j^*)(c)), \\
G_{12}(a\otimes b\otimes c) = \dbl \dbl a,b\dbr ,c \dbr _L
 = \sum\limits_{i,j} e_i (e_j(a))\otimes R(e_j^*)(b) \otimes R(e_i^*)(c).
\end{gather*}
Then for $x,y\in I$ we have
$$
(\langle x,\cdot \rangle \otimes \langle y,\cdot \rangle \otimes \id)F_{23} 
 = R(y)R(x),\quad
(\langle x,\cdot \rangle \otimes \langle y,\cdot \rangle \otimes \id)G_{12} 
 = R(R^*(y)x).
$$

Thus, the identities~\eqref{antiCom},~\eqref{Jacobi} hold if and only if
$R$ is a skew-symmetric RB-operator of weight~0 from $I$ to $\End_f(V)$.
\hfill $\square$

{\bf Remark 3}.
We restrict $R$ in Theorem~2 as an operator from $I$ to $\End_f(V)$ 
instead of $\End(V)$, since otherwise the term $R(yR(x))$ in~\eqref{RB} 
is not well-defined.

{\bf Example 4}~\cite{DoubleLie,DoublePoisson}.
The space $V = F[t]$ equipped with
$$
\dbl t^n,t^m\dbr = \frac{(t^n\otimes 1-1\otimes t^n)
  (t^m\otimes 1-1\otimes t^m)}{t\otimes 1-1\otimes t}
$$
is a double Lie algebra~$L_1$.

Compute the operator $R_1\colon I\to \End(V)$
corresponding to the double Lie algebra~$L_1$,
\begin{equation}\label{R1}
R_1 (e_{ij})
 = \begin{cases}
 - (e_{i,j+1} + e_{i+1,j+2} + \ldots), & i>j \\
 e_{0,j-i+1} + e_{1,j-i+2} + \ldots + e_{i-1,j} , & i\leq j,
 \end{cases}
\end{equation}
where the sum is infinite when $i>j$.
By the formula, $R_1\in\End{}^\prime(V)$ and by Theorem~2, $R_1$
is a skew-symmetric RB-operator from~$I$ to $\End_f(V)$.

Let us identify the matrix algebra $M_n(F)$ of order~$n$ with
$e_{ij}\in I$, $0\leq i,j\leq n-1$.
Given an operator $P$ from $I$ to $\End_f(V)$,
by the projection $P_n$ we mean a linear operator of the space
$\Span\{e_{ij}\mid 0\leq i,j\leq n-1\}$ acting as follows:
$P(e_{ij}) - P_n(e_{ij})\in \Span\{e_{kl}\mid n\leq k \mbox{ or } n\leq l\}$.

One can check that the linear operator $(R_1)_n$ of $M_n(F)$
is an RB-operator of weight~0 on $M_n(F)$ for each~$n$.
Moreover, $(((R_1)_n)^{(\psi_n)})^{(T)}$ coincides with the RB-operator
from~\cite[Example 5.15]{Unital} and it appears 
in~\cite[Example 2.3.3]{Aguiar01}
in terms of the solution of associative Yang---Baxter equation.
Here $\psi_n$ is the automorphism of $M_n(F)$ defined as follows,
$$
\psi(e_{ij}) = e_{n-1-i,n-1-j}.
$$

{\bf Example 5}.
Consider $R_2\in\End{}^\prime(V)$ such that
\begin{equation}\label{R2}
R_2 (e_{ij})
 = \begin{cases}
 -(e_{i-1,j} + e_{i-2,j-1} + \ldots + e_{i-1-j,0}), & i>j, \\
 e_{i,j+1} + e_{i+1,j+2} + \ldots, & i\leq j.
 \end{cases}
\end{equation}
We have defined $R_2$ in such a way that 
$(R_2)_n = (((R_1)_n)^{(\psi_n)})^{(T)}$.
Such definition does not guarantee that we necessarily
obtain a~skew-symmetric RB-operator of weight~0 from~$I$ to $\End_f(V)$.
Thus, we have to state this property of $R_2$.

{\bf Proposition 3}.
The operator $R_2$ is a skew-symmetric RB-operator 
of weight~0 from~$I$ to $\End_f(V)$.

{\sc Proof}.
Firstly, we check the identity
$R_2(e_{ij})R_2(e_{kl}) = R_2( R_2(e_{ij})e_{kl} + e_{ij}R_2(e_{kl}) )$
considering different cases of the values of indices.

{\sc Case 1}: $i>j$, $k>l$. Then
$$
\alpha = R_2(e_{ij})R_2(e_{kl})
 = (e_{i-1,j}+\ldots+e_{i-1-j,0})(e_{k-1,l}+\ldots+e_{k-1-l,0}),
$$

\vspace{-0.9cm}
\begin{multline*}
\beta = R_2(R_2(e_{ij})e_{kl} + e_{ij}R_2(e_{kl})) \\
 = - R_2( (e_{i-1,j}+\ldots+e_{i-1-j,0})e_{kl}
  + e_{ij}(e_{k-1,l}+\ldots+e_{k-1-l,0}) ).
\end{multline*}

Let $k>j$, i.\,e., $j = k-1-p$ for some $p\geq0$. Then
\begin{multline*}
\beta
 = - \chi_{j+l+1-k\geq0}R_2(e_{ij}e_{k-1-p,l-p})
 = - \chi_{j+l+1-k\geq0}R_2(e_{i,j+l+1-k}) \\
 = \chi_{j+l+1-k\geq0}(e_{i-1,j+l+1-k} + \ldots + e_{i-j-l+k-2,0}),
\end{multline*}
since $i>j+l+1-k$.
Here $\chi_{P} = 1$, if $P$ is true, and $\chi_P = 0$, else.
Moreover,
$\alpha = \chi_{j+l+1-k\geq0}(e_{i-1,j+l+1-k} 
 + \ldots + e_{i-j-l+k-2,0}) = \beta$.

Let $k\leq j$, then applying the inequality $i>j+l+1-k$, we compute
$$
\beta
 = - R_2(e_{i-1-j+k,k}e_{kl})
 = - R_2(e_{i-1-j+k,l})
 = e_{i-2-j+k,l} + \ldots + e_{i-2-j+k-l,0}.
$$
On the other hand,
$$
\alpha
 = e_{i-2-j+k,l} + e_{i-2-j+k-1,l-1} + \ldots + e_{i-2-j+k-l,0}
  = \beta.
$$

{\sc Case 2}: $i\leq j$, $k\leq l$. Then
\begin{gather*}
\alpha = R_2(e_{ij})R_2(e_{kl})
 = (e_{i,j+1}+\ldots)(e_{k,l+1}+\ldots), \\
\beta = R_2(R_2(e_{ij})e_{kl} + e_{ij}R_2(e_{kl}))
 = R_2( (e_{i,j+1}+\ldots)e_{kl}
  + e_{ij}(e_{k,l+1}+\ldots) ).
\end{gather*}

Let $k>j$, i.\,e., $j = k-p$ for some $p>0$. Then
$\alpha = e_{i+p-1,l+1} + \ldots = e_{i+k-j-1,l+1} + \ldots$.
Also,
$$
\beta
 = R_2(e_{i+p-1,k}e_{kl})
 = R_2(e_{i+k-j-1,l})
 = e_{i+k-j-1,l+1} + \ldots
 = \alpha,
$$
since $i+k\leq j+l+1$.

Let $j\geq k$, i.\,e., $j = p+k$ for some $p\geq0$. Then
$\alpha = e_{i,l+p+2}+\ldots = e_{i,j+l+2-k}+\ldots$.
Further,
$$
\beta
 = R_2( e_{ij}e_{k+p,l+p+1})
 = R_2( e_{i,l+p+1})
 = R_2( e_{i,l+j-k+1})
 = e_{i,l+j-k+2} + \ldots
 = \alpha,
$$
since $i+k<l+j+2$.

{\sc Case 3}: $i>j$, $k\leq l$. Then
\begin{gather*}
\alpha = R_2(e_{ij})R_2(e_{kl})
 = -(e_{i-1,j}+\ldots+e_{i-1-j,0})(e_{k,l+1}+\ldots), \\
\beta = R_2(R_2(e_{ij})e_{kl} + e_{ij}R_2(e_{kl}))
 = R_2( -(e_{i-1,j}+\ldots+e_{i-1-j,0})e_{kl}
  + e_{ij}(e_{k,l+1}+\ldots) ).
\end{gather*}

When $k>j$, we get $\alpha = \beta = 0$.
Let $j = k+p$ for some $p\geq0$. We compute
$\alpha = -(e_{i-1,l+j-k+1}+\ldots+e_{i-1-j+k,l+1})$.
On the other hand,
$\beta = R_2(-e_{i-1-j+k,l} + e_{i,l+j-k+1})$.
If $i+k\leq j+l+1$, then
$$
\beta
 = -(e_{i-1-j+k,l+1}+\ldots) + (e_{i,l+j-k+2}+\ldots)
 = -(e_{i-1,l+j-k+1}+\ldots+e_{i-1-j+k,l+1})
 = \alpha.
$$
Else,
$$
\beta
 = (e_{i-2-j+k,l}+\ldots+e_{i-2-j+k-l,0})
 - (e_{i-1,l+j-k+1}+\ldots+e_{i-2-j+k-l,0})
 = \alpha.
$$

{\sc Case 4}: $i\leq j$, $k>l$. Then
\begin{gather*}
\alpha = R_2(e_{ij})R_2(e_{kl})
 = -(e_{i,j+1}+\ldots)(e_{k-1,l}+\ldots+e_{k-1-l,0}), \\
\beta = R_2(R_2(e_{ij})e_{kl} + e_{ij}R_2(e_{kl}))
 = R_2( (e_{i,j+1}+\ldots)e_{kl}
  - e_{ij}(e_{k-1,l}+\ldots+e_{k-1-l,0}) ).
\end{gather*}

When $j+1\geq k$, we have $\alpha = 0 = \beta$.

Let $j+2\leq k$, i.\,e., $k = j+2+p$ for some $p\geq0$.
Thus, $\alpha = -(e_{i,l-k+j+2} + \ldots + e_{i+k-j-2,l})$.
Also,
$\beta = R_2(e_{i+k-j-1,l}) - R_2(e_{i,l-k+j+1})$.
If $i+k>j+l+1$, we have
\begin{multline*}
\beta
 = - (e_{i+k-j-2,l}+\ldots+e_{i+k-j-l-2,0})
 + (e_{i-1,l-k+j+1}+\ldots+e_{i+k-j-l-2,0}) \\
 = -(e_{i+k-j-2,l} + \ldots + e_{i,l-k+j+2})
 = \alpha.
\end{multline*}
If $i+k\leq j+l+1$, we have
$$
\beta
 = (e_{i+k-j-1,l+1}+\ldots)
 - (e_{i,l-k+j+2}+\ldots)
 = -(e_{i+k-j-2,l} + \ldots + e_{i,l-k+j+2})
 = \alpha.
$$

Now, we check that $R_2$ is also skew-symmetric operator 
from~$I$ to $\End_f(V)$. Thus, we have to show that
$R_2(e_{ij}) + R_2^*(e_{ij}) = 0$
for all $i,j\geq0$.
By the definition,
$$
R_2^*(e_{ij}) = \sum\limits_{k,l\geq0}R(e_{lk})|_{e_{ji}}e_{kl},
$$
where $R(e_{lk})|_{e_{ji}} = (R(e_{lk}),e_{ij})$
denotes the $e_{ji}$-coordinate of $R(e_{lk})$.

{\sc Case 1}: $i\leq j$. 
Then $R(e_{lk})|_{e_{ji}}$ is nonzero only when
$(l,k)\in\{(j+1+p,i+p)\mid p\geq0\}$.
Thus,
$R_2^*(e_{ij})
 = - (e_{i,j+1}+\ldots)
 = - R_2(e_{ij})$,
as required.

{\sc Case 2}: $i>j$.
Then $R(e_{lk})|_{e_{ji}}$ is nonzero only for
$(l,k)\in\{(j-p,i-1-p)\mid p=0,\ldots,j\}$.
Therefore,
$R_2^*(e_{ij})
 = e_{i-1,j}+\ldots+e_{i-1-j,0}
 = - R_2(e_{ij})$.
\hfill $\square$

{\bf Corollary}.
We have a double Lie algebra structure $L_2$ on~$V$
defined due to~\eqref{Bracket_via_RBInf} by~$R_2$,
$$
\dbl t^n,t^m\dbr 
 = -\frac{(t^n\otimes t^m-t^m\otimes t^n)}{t\otimes 1-1\otimes t}.
$$

Now, we prove that the obtained double Lie algebra~$L_2$ is simple.
It is the first example of a~simple double Lie algebra.

{\bf Theorem 3}.
The double Lie algebra $L_2$ is simple.

{\sc Proof}.
Suppose that $J$ is a~nonzero proper ideal in $L_2$.
Define $n$ as the minimal degree in~$t$ of elements from~$J$.
Let us show that $n = 0$.
If $n>0$, then consider
$f = t^n + \sum\limits_{j=0}^{n-1}\alpha_j t^j\in J$.
We have that the product
$$
\dbl 1,f\dbr
 = t^{n-1}\otimes 1 + t^{n-2}\otimes t + \ldots + 1\otimes t^{n-1}
 + \sum\limits_{j=1}^{n-1}\alpha_j (t^{j-1}\otimes 1 + \ldots + 1\otimes t^{j-1})
$$
lies in $V\otimes J + J\otimes V$.

Consider the map $\psi\colon V\otimes V\to V/J\otimes V/J$
acting as follows: $\psi(v\otimes w) {=} (v+J)\otimes (w+J)$.
On the one hand, $1+J,t+J,\ldots,t^{n-1}+J$ are linearly independent
elements of $V/J$. On the other hand,
$J\otimes V + V\otimes J = \ker(\psi)$.
Thus, $\dbl 1,f\dbr$ is at the same time zero and nonzero
element of $V/J\otimes V/J$.
We obtain a contradiction.
So, $n = 0$ and $1\in J$.

Let us prove by induction on~$s\geq0$ that $t^s\in J$.
For $s = 0$, it is true. Suppose that $s>0$ and we have proved that $t^j\in J$
for all $q<s$. Since $1\in J$, we have
$$
\dbl 1,t^{2s+1}\dbr
 = t^{2s}\otimes 1 + t^{2s-1}\otimes t + \ldots + t^{s+1}\otimes t^{s-1}
 + t^s\otimes t^s + \ldots + 1\otimes t^{2s}\in V\otimes J + J\otimes V.
$$
So, $t^s\otimes t^s\in V\otimes J + J\otimes V$.
Hence, $\psi(t^s\otimes t^s) = 0$, it means that $t^s \in J$.
\hfill $\square$

In the next two examples we consider conjugation of $R_1$ and $R_2$
with transpose and corresponding double Lie algebras.

{\bf Example 6}.
For $R_3 = R_1^{(T)}$, we get a double Lie algebra $L_3$
with the double bracket
$$
\dbl t^n,t^m\dbr
 = \frac{(t^{n+1}\otimes t^{m+1}-t^{m+1}\otimes t^{n+1})}{t\otimes 1-1\otimes t}.
$$

{\bf Example 7}~\cite{DoublePoisson}.
For $R_4 = R_2^{(T)}$, we get a double Lie algebra $L_4$
with the double bracket
$$
\dbl t^n,t^m\dbr = -\frac{(t^{n+1}\otimes1-1\otimes t^{n+1})
 (t^{m+1}\otimes1-1\otimes t^{m+1})}{t\otimes 1-1\otimes t}.
$$

In~\cite{DoublePoisson}, it was stated that each homogeneous double Poisson
algebra on $F[t]$ up to an equivalence is either $L_1$ or $L_4$.
It is easy to show that the double Lie algebras~$L_2$ and~$L_3$
do not satisfy~\eqref{Leibniz}, for example, since
$\dbl t,1\dbr \neq 0$, and therefore do not define the
structure of a double Poisson algebra on $F[t]$.
Note the following connections between double brackets in~$L_1,L_2,L_3,L_4$:
$$
\dbl t^n,t^m\dbr_{L_3} = -(t\otimes t)\dbl t^n,t^m\dbr_{L_2},\quad
\dbl t^n,t^m\dbr_{L_4} = - \dbl t^{n+1},t^{m+1}\dbr_{L_1}.
$$

{\bf Example 8} (V. Kac, see \cite{DoubleLie}).
Consider the double Poisson algebra
$dY(N)= F[t]\otimes M_N(F)$.
Its double bracket relative to the basis
$T_n^{ij} = t^n\otimes e_{ij}$, $n\ge 0$, $i,j=1,\dots, N$,
has the following form:
$$
\dbl T_m^{ij}, T_n^{kl}\dbr
 = \sum\limits_{r=0}^{\min \{m,n\}-1} \big(T_r^{kj}\otimes T^{il}_{m+n-r-1}
 - T^{kj}_{m+n-r-1}\otimes T^{il}_r \big),
$$
the inner bimodule $dY(N)$-action is the associative product.

It is worth mentioning that these relations are similar 
to the defining relations of the Yangian $Y(gl_N)$:
$$
\big[ T_m^{ij}, T_n^{kl} \big]
 = \sum\limits_{r=0}^{\min \{m,n\}-1}
 \big (T_r^{kj} T^{il}_{m+n-r-1} - T^{kj}_{m+n-r-1} T^{il}_r \big ).
$$
We get an RB-operator
$R\colon I\otimes M_N(F)\to \End_f(V)\otimes M_N(F)$ 
such that the double bracket on $dY(N)$
is defined by~\eqref{Bracket_via_RBInf} with the help of~$R$. We have
$$
R(e_{ij}\otimes e_{st})
 = \begin{cases}
 (e_{i,j+1} + e_{i+1,j+2} + \ldots)\otimes e_{st}, & i>j \\
 -(e_{0,j-i+1} + e_{1,j-i} + \ldots + e_{i-1,j})\otimes e_{st} , & i\leq j,
 \end{cases}
$$
where $e_{ij}\in I$, $e_{st}\in M_N(F)$.
Actually, $R = (-R_1)\otimes \id$.

{\bf Remark 4}.
We may extend the double Lie algebra structures $L_1,L_2,L_3,L_4$
on $F[t,t^{-1}]$ and $dY(N)$ on $F[t,t^{-1}]\otimes M_N(F)$ respectively.
It is enough to let both sums in the definition of the corresponding
RB-operator $R_1,R_2,R_3,R_4$, and $R$ be infinite.
Introduce $e_{ij}\in \End(V)$, where
$V = F[t,t^{-1}]$, $i,j\in\mathbb{Z}$, in such a way that
$e_{ij}t^k = \delta_{jk}t^i$. For example, let us extend $R_2$. We define
$$
\widetilde{R}_2(e_{ij})
 = \begin{cases}
 - \sum\limits_{p=0}^{\infty}e_{i-1-p,j-p}, & i>j, \\
   \sum\limits_{p=0}^{\infty}e_{i+p,j+1+p}, & i\leq j.
 \end{cases}
$$
Analogously to the proof of Proposition~3, one
can check that $\widetilde{R}_2$ is a skew-symmetric
RB-operator from $\Span\{e_{ij}\mid i,j\in\mathbb{Z}\}$ to $\End_f(V)$.
By Theorem~2, we get a~double Lie algebra structure $\widetilde{L}_2$ 
on $F[t,t^{-1}]$. Analogously we get double Lie algebras $\widetilde{L}_1$, 
$\widetilde{L}_3$, $\widetilde{L}_4$, and $\widetilde{dY}(N)$.

{\bf Remark 5}.
The operator~$R_2$ is injective. Moreover, $I\subset \Imm(R_2)$.
So, we may define the inverse map $d = R_2^{-1}$ from~$I$ to~$I$.
Then $d(e_{ij}) = e_{i,j-1} - e_{i+1,j}$ is a derivation of~$I$.
By~\cite{DerFinitary}, every such derivation is an inner derivation,
i.\,e., of the form $x\to xa - ax$, where $a\in \End(V)$ such that
$a^{-1}[U]$ is finite-dimensional for each finite-dimensional subspace $U\subset V$.

It is easy to show that $d(x) = Ax - xA$ for $A = e_{10} + e_{21} + \ldots$.
Given a positive integer~$k$, we may define $d_k\colon I\to I$ such that
$d_k(x) = A^k x - xA^k$.
By $d_k$, we may define $P_k\in\End{}^\prime(V)$ as an analogue of
the operator~$R_2$. And all of~$P_k$ provide by~\eqref{Bracket_via_RBInf}
some double Lie algebra structures on $F[t]$.

\section*{Acknowledgements}

The author thanks Korea Institute for Advanced Study (KIAS) 
for the warm hospitality during August 2019,
where a part of the work was done.
The author thanks the anonymous reviewer for the valuable remarks.

The author is supported by the grant of the President 
of the Russian Federation for young scientists (MK-1241.2021.1.1).

\medskip
\noindent Vsevolod Gubarev \\
Sobolev Institute of Mathematics \\
Acad. Koptyug ave. 4, 630090 Novosibirsk, Russia \\
Novosibirsk State University \\
Pirogova str. 2, 630090 Novosibirsk, Russia \\
e-mail: wsewolod89@gmail.com


\begin{thebibliography}{67}
\bibitem{Aguiar00}
M. Aguiar,
Pre-Poisson algebras, Lett. Math. Phys. {\bf 54} (2000) 263--277.

\bibitem{Aguiar00-2}
M. Aguiar,
Infinitesimal Hopf algebras, Contemp. Math. {\bf 267} (2000) 1--30.

\bibitem{Aguiar01}
M. Aguiar,
On the Associative Analog of Lie Bialgebras, J. Algebra {\bf 244} (2001) 492--532.

\bibitem{Finitary}
A. A. Baranov,
Finitary Simple Lie Algebras. J. Algebra {\bf 219} (1999) 299--329.

\bibitem{Baxter}
G. Baxter,
An analytic problem whose solution follows from a simple algebraic identity,
Pacific J. Math. {\bf 10} (1960) 731--742.

\bibitem{BelaDrin82}
A.~A. Belavin, V.~G. Drinfel'd,
Solutions of the classical Yang---Baxter equation for simple Lie algebras,
Funct. Anal. Appl. (3) {\bf 16} (1982) 159--180.

\bibitem{Crawley-Boevey}
W. Crawley-Boevey,
Poisson structures on moduli spaces of representations,
J.~Algebra (1) {\bf 325} (2011) 205--215.

\bibitem{DoubleLie}
M. E. Goncharov, P. S. Kolesnikov.
Simple finite-dimensional double algebras, J.~Algebra {\bf 500} (2018) 425--438.

\bibitem{Unital}
V. Gubarev,
Rota---Baxter operators on unital algebras. Moscow Math. Journal (accepted),
arXiv.1805.00723v3.

\bibitem{GuoMonograph}
L. Guo,
An Introduction to Rota---Baxter Algebra. Surveys of Modern Mathematics, vol. 4,
International Press, Somerville (MA, USA); Higher education press, Beijing, 2012.

\bibitem{RBModules}
L. Guo and Z. Lin,
Representations and modules of Rota---Baxter algebras, preprint (2015),
arXiv:1905.01531, 28~p.

\bibitem{IKV}
N. Iyudu, M. Kontsevich, Y. Vlassopoulos,
Pre-Calabi-Yau algebras as noncommutative Poisson structures,
J. Algebra {\bf 567} (2021) 63--90.

\bibitem{DerFinitary}
S.~G. Kolesnikov, N.~V. Mal'tsev,
Derivations of matrix rings containing a subring of triangular matrices,
Russian Math. (Iz. VUZ) {\bf 55} (11) (2011) 18--26.

\bibitem{O-operator}
B.~A.~Kupershmidt,
What a Classical $r$-Matrix Really is,
J. Nonlinear Math. Phy. {\bf 6} (1999) 448--488.

\bibitem{DoublePoissonFree}
A. Odesskii, V. Rubtsov, V. Sokolov,
Double Poisson brackets on free associative algebras,
in: Noncommutative birational geometry, representations and combinatorics,
Contemp. Math. {\bf 592} (2013) 225--239, Amer. Math. Soc., Providence, RI.

\bibitem{PichereauWeyer}
A. Pichereau and G. Van den Weyer,
Double Poisson cohomology of path algebras of quivers, 
J.~Algebra {\bf 319} (2008) 2166--2208.

\bibitem{Polishchuk}
A. Polishchuk,
Classical Yang-Baxter equation and the $A_\infty$-constraint,
Adv. Math. (1) {\bf 168} (2002) 56--95.

\bibitem{Schedler}
T. Schedler,
Poisson algebras and Yang-Baxter equations, in:
Advances in Quantum Computation (Contemp. Math., Vol. 482),
AMS, Providence, R.I. (2009) 91--106.

\bibitem{Sokolov}
V.~V. Sokolov,
Classification of constant solutions of the associative
Yang---Baxter equation on $\mathrm{Mat}_3$,
Theor. Math. Phys. (3) {\bf 176} (2013) 1156--1162.

\bibitem{Kac15}
A. De Sole, V. G. Kac, D. Valeri,
Double Poisson vertex algebras and non-commutative Hamiltonian equations,
Adv. Math. {\bf 281} (2015) 1025--1099.

\bibitem{Tricomi}
F.~G.~Tricomi,
On the finite Hilbert transformation.
Quart. J. Math. {\bf 2} (1951) 199--211.

\bibitem{Uchino}
K.~Uchino,
Quantum Analogy of Poisson Geometry, Related Dendriform Algebras
and Rota-Baxter Operators, Lett. Math. Phys {\bf 85} (2008) 91--109.

\bibitem{DoublePoisson}
M. Van den Bergh,
Double Poisson algebras, Trans. Amer. Math. Soc. {\bf 360} (11) (2008),5711--5769.

\bibitem{Zhelyabin}
V.N. Zhelyabin,
Jordan bialgebras of symmetric elements and Lie bialgebras,
Siberian Mat. J. (2) {\bf 39} (1998) 261--276.
\end{thebibliography}
\end{document}